\documentclass[a4paper,DIV=12]{scrartcl}
\usepackage[utf8x]{inputenc}
\usepackage{amsmath,amssymb,graphicx,subfig}
\usepackage[english]{babel}

\usepackage{hyperref}
\hypersetup{pdfauthor={W. Huang, L. Kamenski, and J. Lang},
   pdftitle={Adaptive finite elements with anisotropic meshes}}

\providecommand{\enorm}[1]{||| #1 |||}
\providecommand{\norm}[1]{\lVert #1 \rVert}

\begin{document}

\title{Adaptive finite elements \\with anisotropic meshes\thanks{Supported in part by the National Science Foundation (U.S.A.) through grant DMS-1115118 and the German Research Foundation through grants SFB568/3, SPP1276 (MetStroem) and KA~3215/1-1.
}}
\author{ 
   Weizhang Huang\thanks{Department of Mathematics, the University of Kansas
      (\href{mailto:huang@math.ku.edu}{huang@math.ku.edu})}
   \and Lennard Kamenski\thanks{Department of Mathematics, the University of Kansas
      (\href{mailto:lkamenski@math.ku.edu}{lkamenski@math.ku.edu})}
   \and Jens Lang\thanks{Department of Mathematics,  Technische Universit{\"a}t Darmstadt
      (\href{mailto:lang@mathematik.tu-darmstadt.de}{lang@mathematik.tu-darmstadt.de})\newline
      Center of Smart Interfaces, Technische Universit{\"a}t Darmstadt\newline
      Graduate School of Computational Engineering, Technische Universit{\"a}t Darmstadt}
}
\date{}
\subject{}
\maketitle

\begin{abstract}
The paper presents a numerical study for the finite element method with anisotropic meshes.
We compare the accuracy of the numerical solutions on quasi-uniform, isotropic, and anisotropic meshes for a test problem which combines several difficulties of a corner singularity, a peak, a boundary layer, and a wavefront.
Numerical experiment clearly shows the advantage of anisotropic mesh adaptation.
The conditioning of the resulting linear equation system is addressed as well.
In particular, it is shown that the conditioning with adaptive anisotropic meshes is not as bad as generally assumed.
\end{abstract}

\section{Introduction}
\label{sec:introduction}

Anisotropic mesh adaptation, i.e., adaptation of the size and shape of mesh elements, has been shown to be of significant advantage for problems with distinct anisotropic features.
Moreover, the ability to adjust the shape and orientation of mesh elements has proven to be useful for designing numerical schemes with particular features, e.g., satisfying the discrete maximum principle \cite{LiHua10} or improving the conditioning of the finite element equations \cite{D'RoDo97,KaHuXu12}.

In this paper we concentrate on obtaining anisotropic meshes for the purpose of minimizing the numerical solution error.
Typically, the optimal shape and orientation of mesh elements depend on the Hessian \cite{D'Azev91,ForPer03,Huang06,Simpso94} or the first derivatives \cite{Picass03,Picass06} of the exact solution of the underlying problem.
This is the first major difficulty since the exact solution is usually not available. 
One possibility to solve this difficulty is to try to recover the approximate Hessian from the numerical solution in the course of the computation.
In \cite{FoMiPe04,ForPer03}, mesh adaptation is based on a residual-based error estimator but still requires Hessian recovery for the solution of the dual problem.
Unfortunately, Hessian recovery methods do work very well for interpolation problems but they cannot provide a convergent recovery if applied to linear finite element approximations on non-uniform meshes \cite{AgLiVa10,Kamens09PhD}, although adaptive finite element methods based on Hessian recovery still provide excellent mesh adaptation in practice \cite{Dolejs98,FoMiPe04,ForPer03,HuaLi10,LiHua10,VasLip99}.
The fact that the convergence of adaptive algorithms employing Hessian recovery cannot be proven directly (since the recovered Hessian does not converge to the exact Hessian) explains recent interest in anisotropic mesh adaptation based on some kind of a~posteriori error estimates which do not depend on the exact solution of the underlying problem \cite{AgLiVa08,AgLiVa09,AgLiVa10,ApGrJM04,CaHuRu01,FrLaRo96,GeHaHo07,HuKaLa10}.
Moreover, as shown in \cite[Sect.~5.3]{HuKaLa10}, using error estimates could be of advantage for problems exhibiting gradient jumps or similar discontinuities along internal interfaces because methods based on recovery of derivatives could result in unnecessarily high mesh density near discontinuities.

For our study we employ the anisotropic mesh adaptation algorithm from \cite{HuKaLa10} which employs a globally defined hierarchical basis error estimate (HBEE) for obtaining the directional information.
In contrast to the recovery-based algorithms, this method adapts the mesh in order to directly minimize the a~posteriori error estimate and, thus, relies on the accuracy of the error estimator but does not require recovery of derivatives of the exact solution.
In this sense, the algorithm is completely a~posteriori.

Another major concern when using anisotropic meshes is the conditioning of the finite element equations.
Generally speaking, an anisotropic mesh is expected to contain elements of large aspect ratio\footnote{In this paper the aspect ratio of a triangular element is defined as the longest edge divided by the shortest altitude. For example, an equilateral triangle has an aspect ratio of $2 / \sqrt{3} \approx 1.15$.} and there exists a concern that an anisotropic mesh will lead to extremely ill-conditioned linear algebraic systems and this may weaken the accuracy improvements gained through anisotropic mesh adaptation.
Fortunately, as it has been recently shown in \cite{KaHuXu12}, the conditioning of the stiffness matrix with  anisotropic meshes is not necessarily as bad as generally assumed, especially in $2D$.
In Sect.~\ref{subsec:conditioning}, we will see that even if the condition number of the stiffness matrix with an anisotropic mesh is larger than that with an isotropic mesh, the accuracy gained through anisotropic mesh adaptation still clearly outbalances the conditioning issues, at least for the example considered.

The outline of this paper is as follows: a brief description of the adaptation algorithm is given in Sect.~\ref{sec:discretization} which is followed by the numerical experiment in Sect.~\ref{sec:numericalExperiment}.
The concluding remarks are given in Sect.~\ref{sec:summary}.

\section{Discretization and the mesh adaptation algorithm}
\label{sec:discretization}

We consider a Dirichlet problem for the Poisson equation
\begin{equation}
   \begin{cases}
      - \Delta u = f, & \text{in } \Omega \\
      u = 0,                            & \text{on } \partial\Omega
   \end{cases}
   \label{eq:bvp}
\end{equation}
where $\Omega \subset \mathbb{R}^2$ is a connected bounded polygonal domain. 

For a given triangulation $\mathcal{T}_h$ of $\Omega$ and the associated linear finite element space $V_h \subset H_0^1(\Omega)$, the linear finite element solution $u_h \in V_h$ of \eqref{eq:bvp} is defined by
\begin{equation}
   \int_{\Omega} \nabla v_h  \cdot \nabla u_h ~ dx 
      = \int_{\Omega} f  v_h ~ dx, \quad \forall v_h \in V_h.
    \label{eq:femIntegral}
\end{equation}
The finite element space $V_h$ and the finite element solution $u_h$ can be written as
\begin{equation}
   V_h = \text{span} \{ \phi_1, \cdots, \phi_{n_{int}} \}
   \quad \text{and} \quad
   u_h = \sum_{j=1}^{n_{int}} u_j \phi_j ,
   \label{eq:uh}
\end{equation}
where $\phi_j$ is the standard linear basis function associated with the $j$-th vertex and $n_{int}$ is the number of interior vertices of the triangulation.
Substituting \eqref{eq:uh} into \eqref{eq:femIntegral} and taking $v_h = \phi_i$ for $i = 1, \dots, n_{int}$ results in the linear system
\begin{equation}
   A u_h = F,
   \label{eq:linearSystem}
\end{equation}
where
\[ A_{ij} =  \int_{\Omega} \nabla \phi_j \cdot \nabla \phi_i ~ dx
   \quad \text{and} \quad
   F_i    =  \int_{\Omega} f  \phi_i ~ dx.
\]
Note that in order to obtain the finite element solution $u_h$ we need to solve the linear algebraic system \eqref{eq:linearSystem}.
Thus, the accuracy of $u_h$ depends also on the conditioning of this system which in turn is affected by the choice of the mesh.
As mentioned in the introduction, there is a concern that anisotropic meshes could lead to extremely ill-conditioned linear systems and this may weaken the accuracy gained with anisotropic mesh adaptation.
In our numerical experiment in Sect.~\ref{subsec:conditioning} we will address this issue in detail.

In order to construct $\mathcal{T}_h$ (and, thus, the corresponding $V_h$) we employ the $M$-uniform mesh approach which generates an adaptive mesh as a quasi-uniform one in the metric specified by a symmetric and strictly positive definite tensor $M = M(x)$ \cite{Huang05a}. 
The algorithm starts with an initial mesh.
For every mesh $\mathcal{T}_h^{(i)}$ we compute the finite element solution $u_h^{(i)}$ which is used to compute a new adaptive mesh for the next iteration step.
The new mesh is generated as an $M$-uniform mesh with a metric tensor $M_h^{(i)}$ computed from $u_h^{(i)}$.
This yields the sequence
\[
   \mathcal{T}_h^{(0)} \rightarrow u_h^{(0)} \rightarrow M_h^{(0)} \rightarrow 
   \mathcal{T}_h^{(1)} \rightarrow u_h^{(1)} \rightarrow M_h^{(1)} \rightarrow \dots 
\]
The mesh adaptation process is repeated until the mesh is $M$-uniform within a given tolerance (see \cite[Sect.~4.1]{HuKaLa10} for more details).
In our computation we use BAMG (\emph{bidimensional anisotropic mesh generator} \cite{bamg}) to construct anisotropic meshes for a given metric tensor $M$.

Typically, the optimal metric tensor $M_h$ depends on the Hessian of the exact solution \cite{ForPer01,ForPer03,Huang05a} which is usually unknown.
In this study we follow \cite{HuKaLa10} and employ the hierarchical basis a~posteriori error estimate (HBEE) to obtain the directional information required for the metric tensor $M_h$.
The brief idea is as follows (see \cite{HuKaLa10} for details).

If we have an error estimate $z_h$ such that
\[
   \norm{u - u_h} \leq C \norm{z_h}.
\]
for a given norm $\norm{\cdot}$ and if it further has the property $\Pi_h z_h \equiv 0$ with $\Pi_h$ being the interpolation operator associated with $V_h$ (which is fulfilled by the HBEE), than the finite element approximation error is bounded by the interpolation error of the error estimate, 
\begin{equation}
  \|u - u_h\| \le C \norm{z_h} = C \norm{z_h - \Pi_h z_h}.
   \label{eq:zhIE}
\end{equation}
Hence, up to a constant, the solution error is bounded by the interpolation error of the error estimate and the mesh can be constructed to minimize the interpolation error of $z_h$; the metric tensor $M_h$ does not depend on the Hessian of the exact solution.

In this study, we are concerned with the error measured in the $H^1$ semi-norm, which is the energy norm from \eqref{eq:bvp}.
Therefore, instead of using the metric tensor developed in \cite{HuKaLa10} for the error measured in the $L^2$ norm, we construct the metric tensor which minimizes the interpolation error of $z_h$ measured in the $H^1$ semi-norm.
In two dimensions the optimal metric tensor is given element-wise by
\[
   M_{K} = \left\|I + \frac{1}{\alpha_h} |H_K(z_{h})|\right\| 
      \cdot \det\left(I + \frac{1}{\alpha_h} |H_K(z_{h})|\right)^{-\frac{1}{4}}
      \cdot \left[I + \frac{1}{\alpha_h} |H_K(z_{h})|\right],
\]
where $H_K(z_h)$ denotes the Hessian of the (quadratic) hierarchical basis error estimate $z_h$ on element $K$ and $\alpha_h$ is a regularization parameter to ensure that $M_K$ is strictly positive definite.
$\alpha_h$ can also be seen as an adaptation intensity control: uniform mesh has $\alpha_h = \infty$ and if $\alpha_h \rightarrow 0$ the mesh becomes more adaptive.
Usually, $\alpha_h$ is chosen so that about half of the mesh elements are concentrated in regions where $\det(M)$ is large (see \cite{Huang06} for more details on the choice of $M_K$ and $\alpha_h$). 

In our computations we employ the globally defined hierarchical basis error estimate since it contains more directional information of the solution than localized versions \cite[Sect.~5.1]{HuKaLa10}.
Moreover, it has been shown that local error estimates can be inaccurate on anisotropic meshes \cite{DoGrPf99}.
To avoid the cost of the exact solution of the global error problem, we use only a few sweeps of the symmetric Gauss-Seidel iteration for the resulting linear system until the relative difference of the old and the new error approximations is under a given relative tolerance.
This proves to be adequate for the purpose of mesh adaptation and the computational cost is comparable to that of the Hessian recovery: in the tests, the computation of HBEE is about twice slower than Hessian recovery.

Although the validity of the classical hierarchical basis error estimate $z_h$ for the anisotropic case is still unclear, theoretical considerations in \cite[Sect.~6.4]{Grosma06PhD} and numerical results in \cite{HuKaLa10} suggest that the hierarchical basis error estimate is a reliable source of information when a mesh is aligned with the solution.

\section{Numerical experiment}
\label{sec:numericalExperiment}

\begin{figure}[t] \centering
   \includegraphics[height=0.20\textheight,clip]{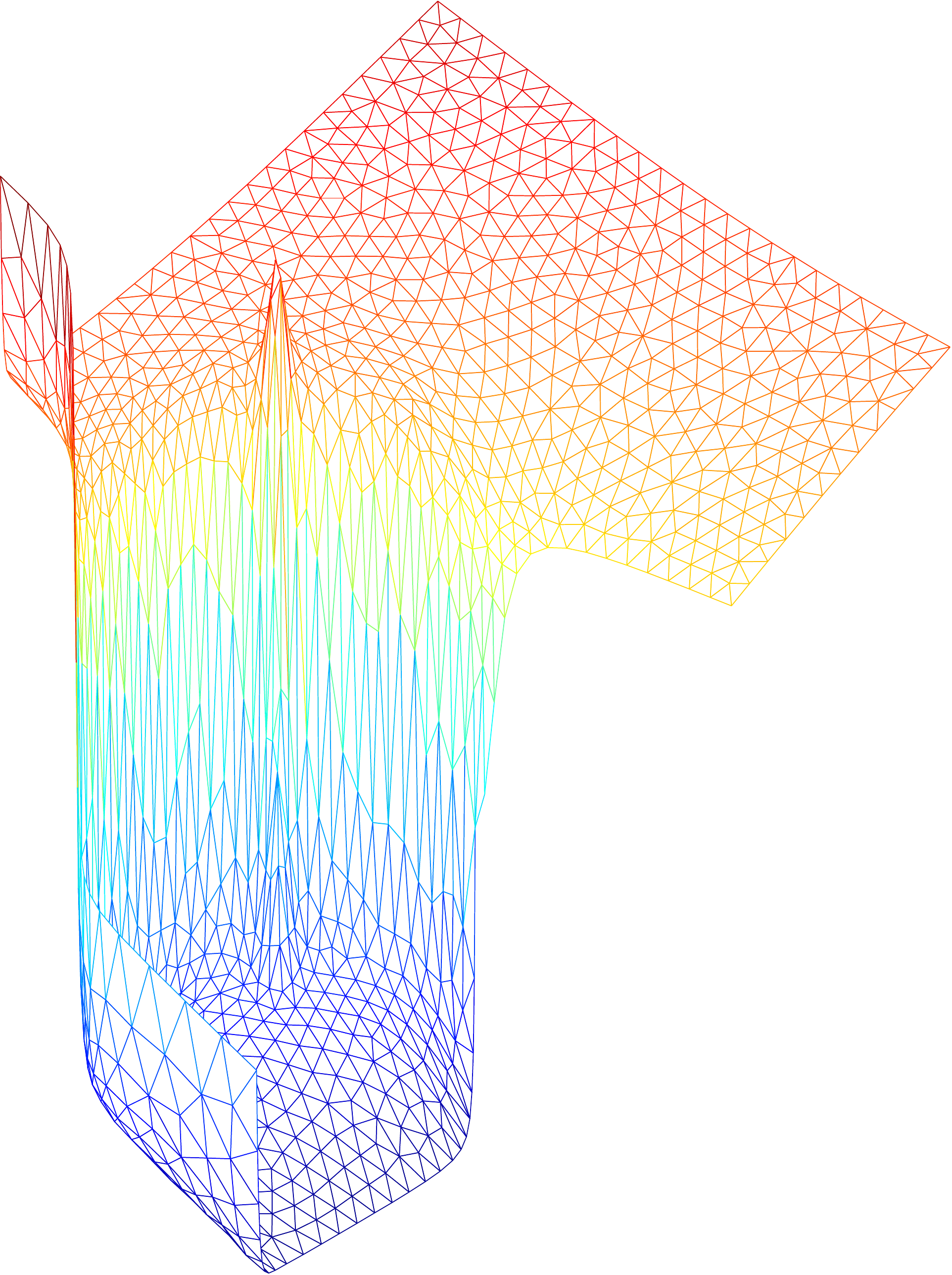}
   \qquad
   \includegraphics[height=0.19\textheight,clip]{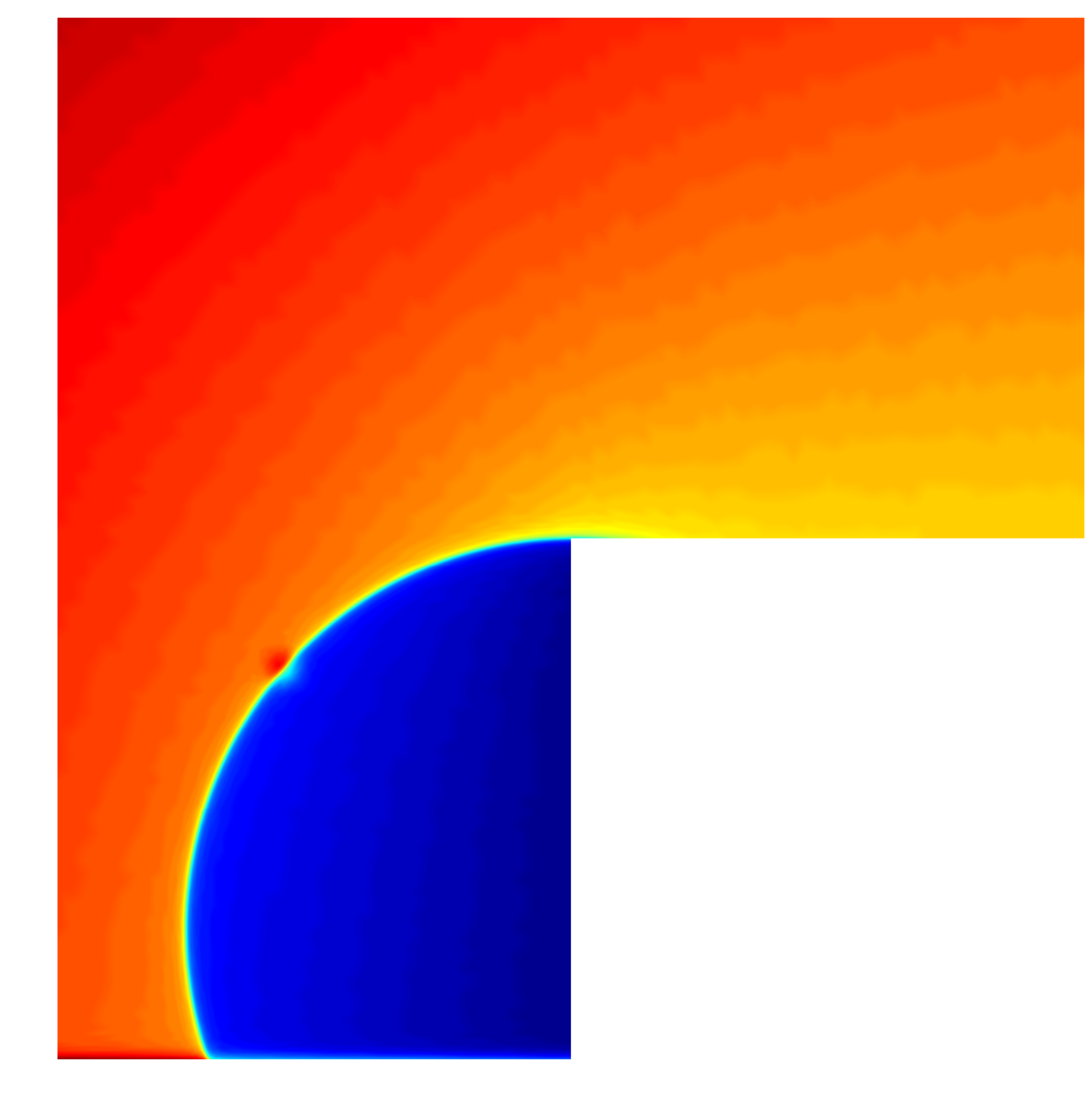}
   \caption{Surface and color plots of the numerical solution.}
   \label{fig:temperaturePlot}
\end{figure}

For the numerical experiment we consider a problem in \cite{Mitche10} which combines multiple difficulties.
It is a Dirichlet problem of the Poisson equation
\begin{equation}
   \begin{cases}
      - \Delta u = f,   & \text{in } \Omega\\
      u = g, & \text{on } \partial\Omega
   \end{cases}
   \label{eq:Poisson}
\end{equation}
where $\Omega$ is an L-shaped domain $\Omega = (-1, 1) \times (-1, 1) \setminus [0, 1) \times (-1, 0]$. 
The functions $f$ and $g$ are chosen such that the exact solution $u$ is given by
\begin{align*}
   u(x,y) 
      &= r^{2/3} \sin(2 \theta /3) + \tan^{-1} \left( 200 \left(\sqrt{x^2 + (y + 3/4)^2} - 3/4 \right)\right) \\
      & \quad + e^{-1000 \left( (x + \sqrt{5}/4)^2 + (y + 1/4)^2 \right)} + e^{- 100 (y + 1)},
\end{align*}
where $r$ and $\theta$ are the polar coordinates.
The solution has 
\begin{itemize}
   \item a singular gradient at $(0,0)$ due to a reentrant corner of the L-shaped domain $\Omega$,
   \item a circular wavefront with the center in $(0,-3/4)$ and the radius of $3/4$,
   \item a sharp peak at $(-\sqrt{5}/4, -1/4)$,
   \item and a boundary layer along the line  $y = -1$.
\end{itemize}
Figure~\ref{fig:temperaturePlot} shows the surface and the color plot of a numerical solution.

\subsection{Accuracy of the numerical solution}
\label{subsec:optimal:mesh:grading}

\begin{figure}[p] \centering
   \subfloat[Delaunay: 2326 elements; max. aspect ratio $2.8$.]{
         \includegraphics[width=0.40\textwidth,clip]{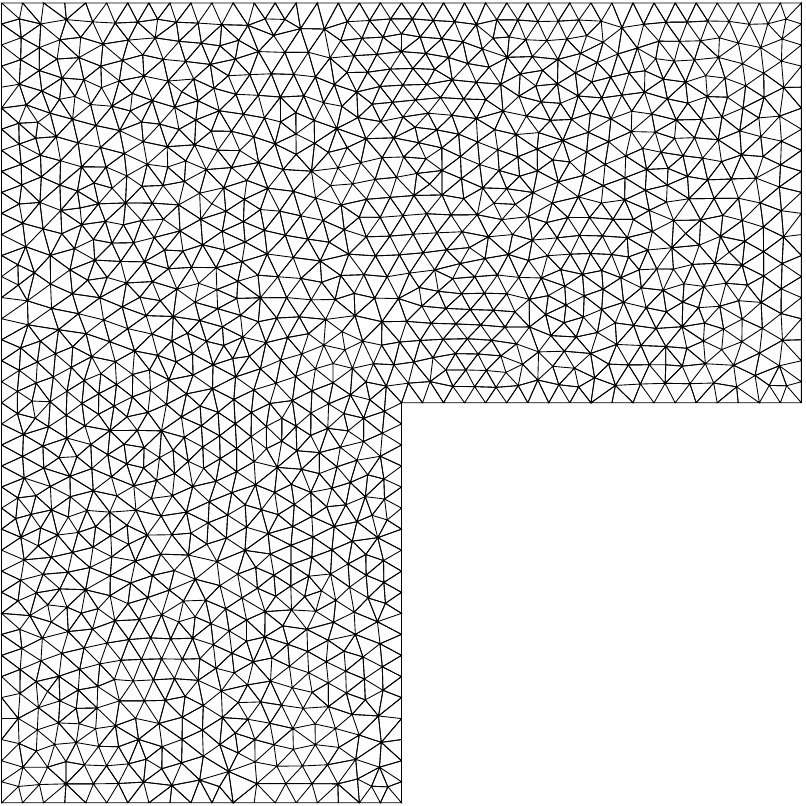}
         \quad
         \includegraphics[width=0.40\textwidth,clip]{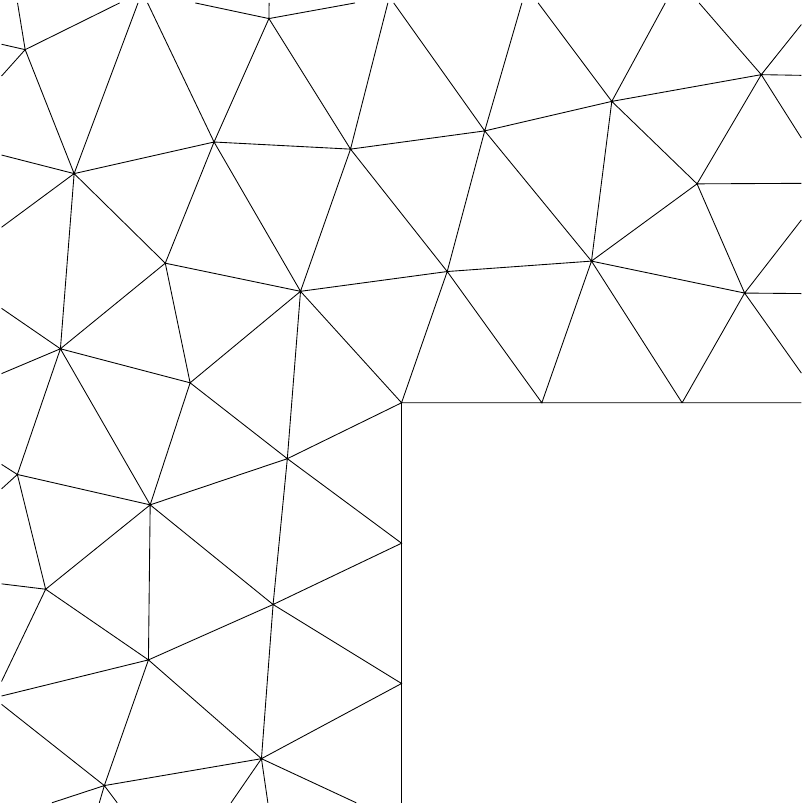}
      \label{fig:meshExamples:Delaunay}
   }\,
    \subfloat[Isotropic adaptive: 2321 elements; max. aspect ratio $3.0$.]{
         \includegraphics[width=0.40\textwidth,clip]{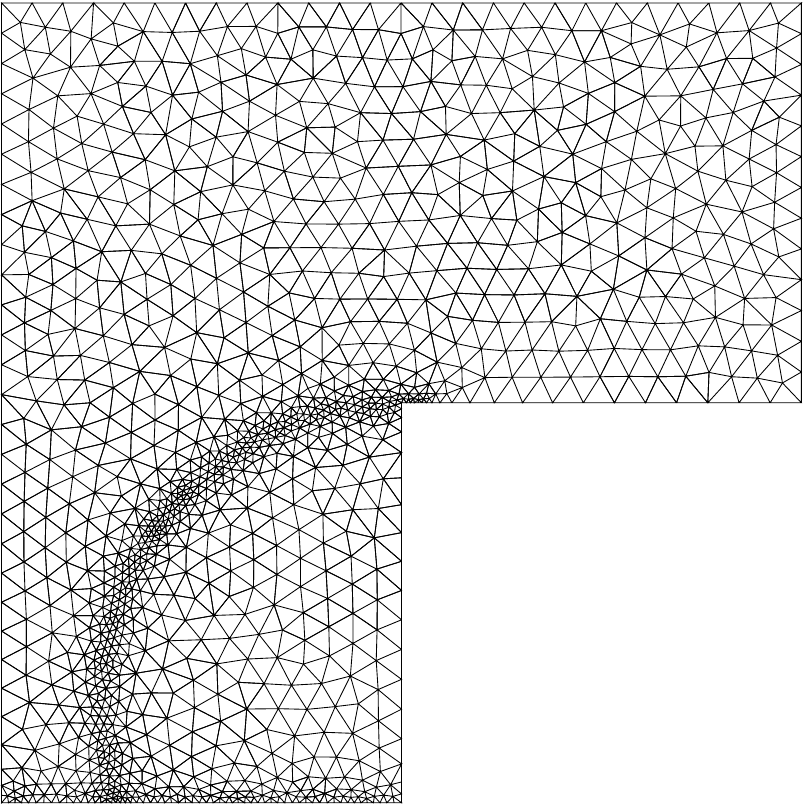}
         \quad
         \includegraphics[width=0.40\textwidth,clip]{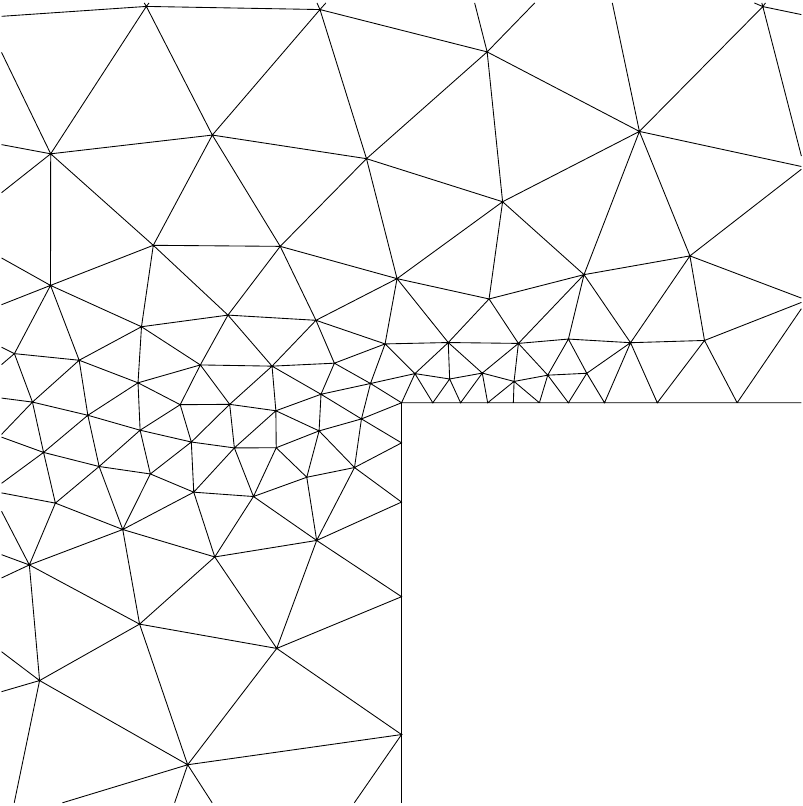}
      \label{fig:meshExamples:isotropic}
   }\,
   \subfloat[Anisotropic adaptive: 2316 elements; max. aspect ratio $24.4$.]{
         \includegraphics[width=0.40\textwidth,clip]{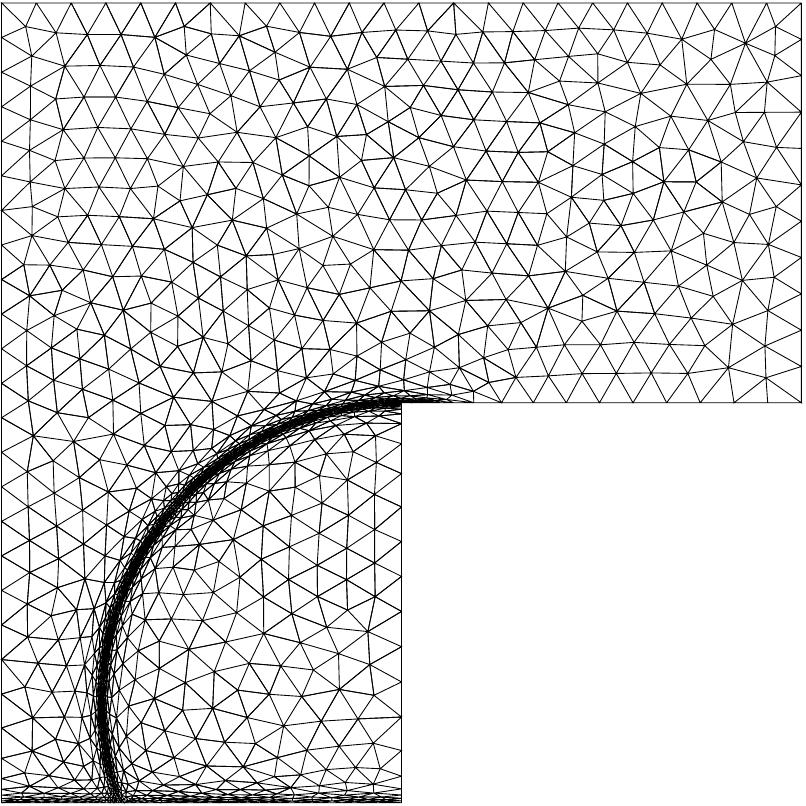}
         \quad
         \includegraphics[width=0.40\textwidth,clip]{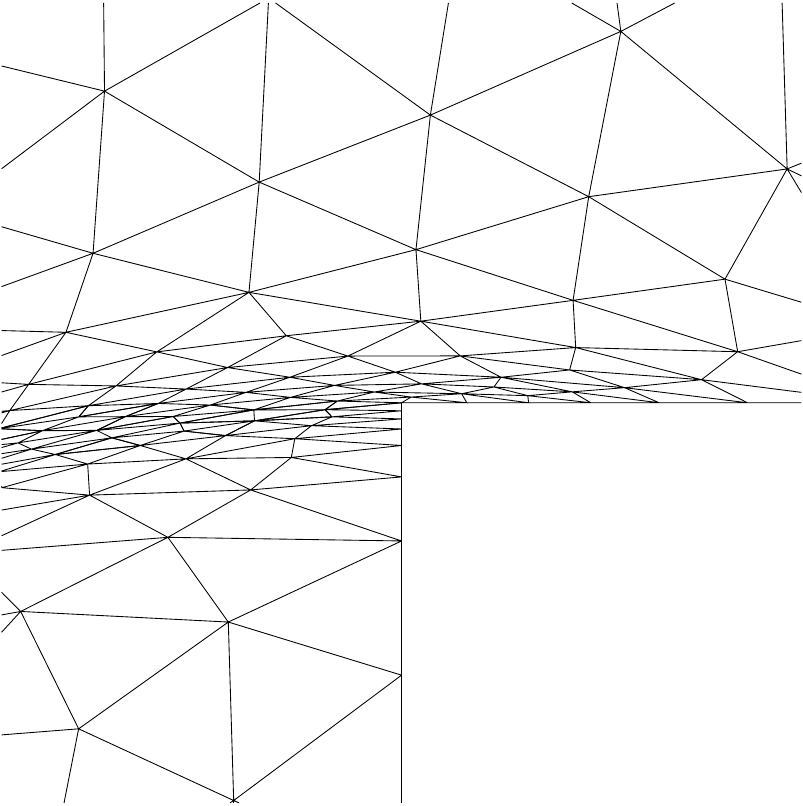}
      \label{fig:meshExamples:anisotropic}
   }
   \caption{Mesh examples and 6.6 times close-up views at the reentrant corner.}
   \label{fig:meshExamples}
\end{figure}

First, we compare the accuracy of the numerical solution for Delaunay (quasi-uniform), adaptive isotropic, and adaptive anisotropic meshes.
Examples of mesh types are given in Fig.~\ref{fig:meshExamples}.
We observe that both isotropic and anisotropic adaptive meshes (Figs.~\ref{fig:meshExamples:isotropic} and \ref{fig:meshExamples:anisotropic}, respectively) have high mesh density in regions with difficulties but the anisotropic mesh (Fig.~\ref{fig:meshExamples:anisotropic}) is clearly much better aligned with the steep boundary layer and the wavefront.
This is the major difference between the isotropic and anisotropic adaptation: the isotropic adaptation can provide proper mesh density whereas the anisotropic adaptation can provide both proper mesh density \emph{and} proper alignment of the mesh with the anisotropic features of the solution.

\begin{figure}[t] \centering
   \includegraphics[width=0.75\textwidth,clip]{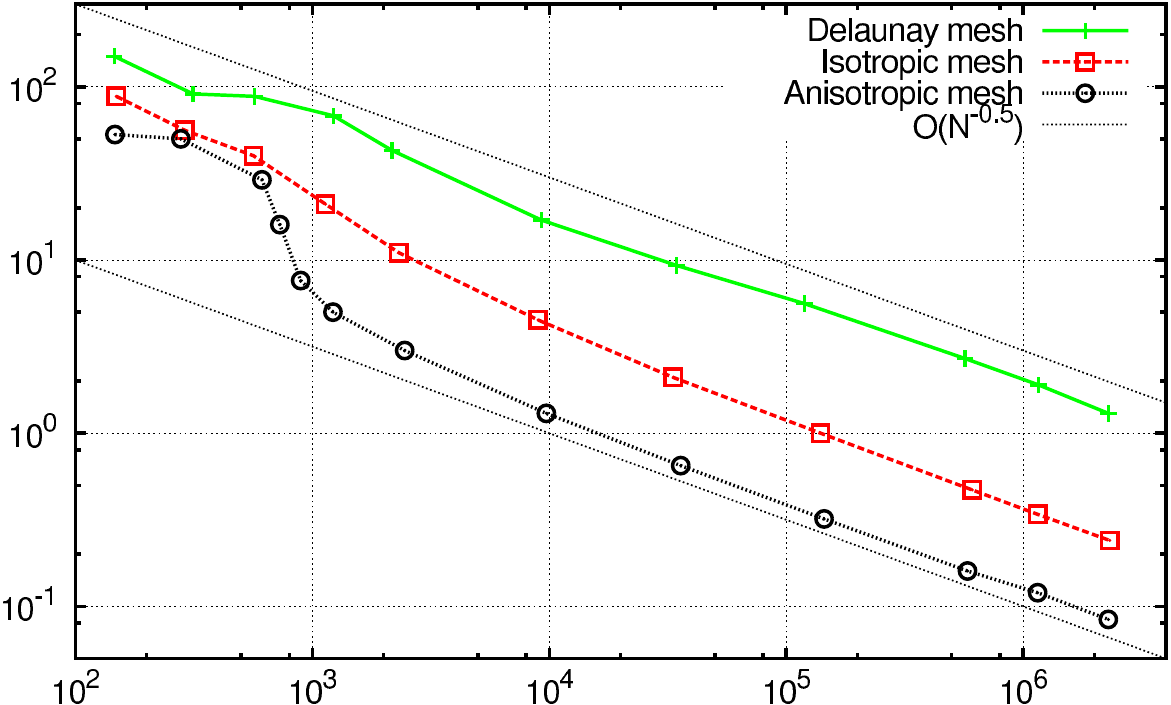}
   \caption{Energy norm of the finite element error vs. number of mesh elements.}
   \label{fig:convergence}
\end{figure}

Fig.~\ref{fig:convergence} shows the error of the numerical solution measured in the energy norm $\enorm{u - u_h}$, which is equal to $H^1(\Omega)$ semi-norm $|u - u_h|_{H^1(\Omega)}$ for the example considered.
The convergence plot shows that an anisotropic adaptive mesh requires ca. 200 times fewer elements than a quasi-uniform mesh in order to achieve the same accuracy and ca. 10 times fewer elements than an isotropic adaptive mesh.
In other words, the finite element solution with an anisotropic mesh has a 15 times smaller error than an error of the solution on a quasi-uniform mesh with the same number of elements and 3 times smaller than the error achieved by means of an isotropic adaptive mesh.
The asymptotic convergence order of the error in the energy norm is the same for all three kinds of meshes: it is $O(N^{-0.5})$.\footnote{Note that $O(N^{-0.5}) = O(h)$ for quasi-uniform meshes in 2D.}
This is expected since we cannot have a better convergence order for anisotropic mesh adaptation but can expect a much smaller constant when the solution of the problem has anisotropic features.
In our test example we gain more than one order of magnitude in comparison to quasi-uniform meshes and about one half of the order in comparison to the isotropic adaptation.

Fig.~\ref{fig:convergence} provides also an interesting insight into the behaviour of anisotropic mesh adaptation.
For very coarse meshes ($N < 300$) the resolution is not good enough to capture the anisotropy of the solution, the mesh is isotropic and has the same error as with isotropic mesh adaptation.
The interesting part of the plot is between $N \approx 300$ and  $N \approx 1000$, where the algorithm starts to catch the anisotropic features and the error drops quickly.
When the anisotropic mesh is fine enough to resolve the anisotropy of the solution ($N > 1000$), the error convergence rate reaches the asymptotic state.

\subsection{Condition number of the stiffness matrix}
\label{subsec:conditioning}
In this section, we compute the exact condition number (with respect to the $\lVert\cdot\rVert_2$ matrix norm) for the stiffness matrix of the anisotropic finite elements equations and compare it to the conditioning of the finite element equations with isotropic adaptive and quasi-uniform meshes.

\begin{figure}[t] \centering
   \subfloat[Unscaled.]{
      \includegraphics[width=0.48\textwidth,clip]{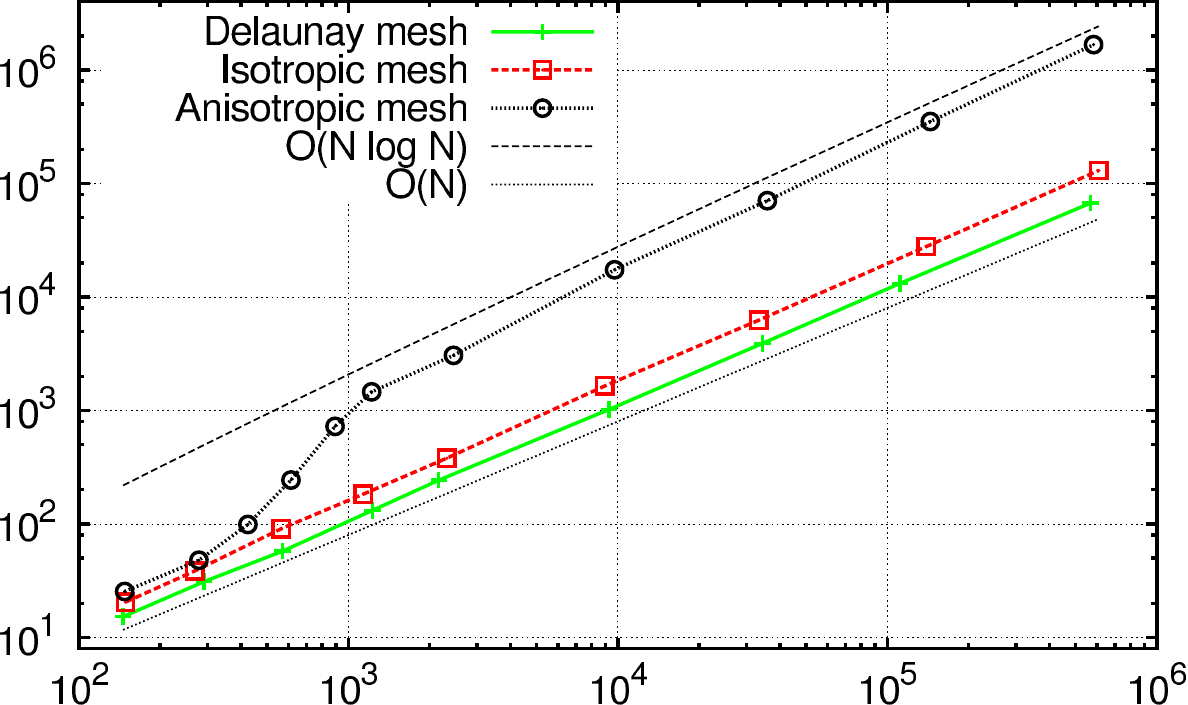}
      \label{fig:condition:unscaled}
   }
   \subfloat[After diagonal scaling.]{
      \includegraphics[width=0.48\textwidth,clip]{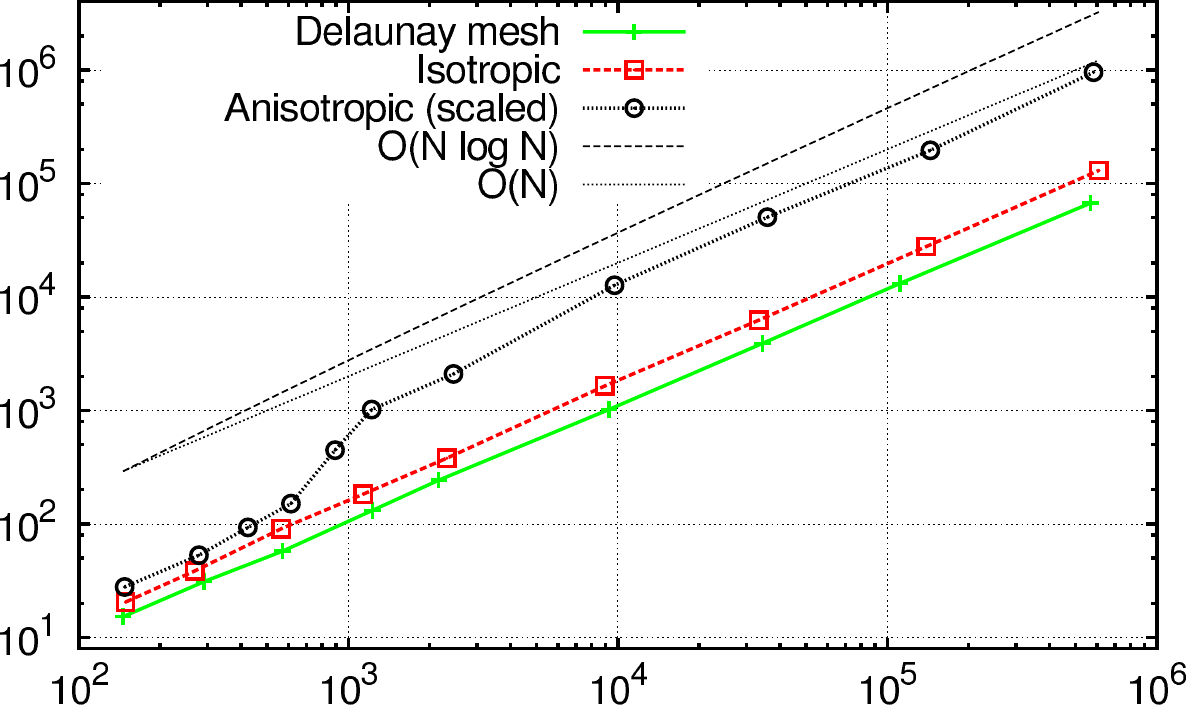}
      \label{fig:condition:scaled}
   }
   \caption{Condition number of the stiffness matrix  vs. number of elements.}
   \label{fig:condition}
\end{figure}

The analysis in \cite{KaHuXu12} for the Laplace operator in 2D shows that the condition number can be bounded by a term depending mainly on the number $N$ and the largest aspect ratio of the mesh elements.
In our numerical experiment the maximum aspect ratio is up to $3.8$ for quasi-uniform and isotropic meshes and up to $37.9$ for anisotropic meshes.
Thus, the rough estimate on the ratio between the condition numbers of the anisotropic and isotropic systems  should be about $37.9/3.8 \approx 10.0 $.
This is in perfect agreement with our numerical results presented in Fig.~\ref{fig:condition:unscaled} which show that the condition number of the linear system with anisotropic meshes is about one order of magnitude higher than that with the isotropic meshes.
Notice also the sudden jump in the condition number for the anisotropic case in the range  $300 \leq N \leq 1000$: the algorithm starts to catch the anisotropic features of the solution and the maximum aspect ratio of the mesh increases quickly as the mesh becomes more and more anisotropic (cp. the corresponding error decrease in Fig.~\ref{fig:convergence}).

Fig.~\ref{fig:condition:unscaled} also shows that that the asymptotic behaviour of the condition number with anisotropic meshes is at most $O(N \log N)$ which is only slightly larger than $O(N)$ in the quasi-uniform case.
The conditioning with isotropic adaptive meshes is also slightly larger than $O(N)$ although still smaller than $O(N \log N)$.
Moreover, if a mesh is only locally anisotropic (as in our example), a proper diagonal scaling can reduce the conditioning of the stiffness matrix so that it is comparable with the condition number in the uniform case (see \cite{KaHuXu12} for more details on diagonal scaling).
Fig.~\ref{fig:condition:scaled} shows that the asymptotic rate of the conditioning of the scaled stiffness matrix is reduced to essentially $O(N)$, which is comparable to that with uniform meshes.

\section{Conclusion}
\label{sec:summary}

Our numerical experiment shows that for problems with anisotropy the anisotropic mesh adaptation is clearly superior to the isotropic one.
In our example, at least a half order of magnitude could be gained in accuracy by switching from the isotropic mesh adaptation to the anisotropic one.
The globally defined hierarchical basis error estimate provides good directional information for the anisotropic mesh generation, provided the number of mesh elements is large enough to resolve the anisotropy of the solution.
It is worth pointing out that the results in Fig.~\ref{fig:convergence} present the error of the final numerical solution, i.e., \emph{after} solving the linear system.
Thus, even if the condition number of the linear system with anisotropic meshes is larger than that with isotropic meshes, the accuracy gained through the anisotropic discretization for problems with anisotropic features outbalances possible losses due to the numerical accuracy.

\subsubsection*{Acknowledgement} 
The authors are thankful to the anonymous referee for the valuable comments.

\bibliographystyle{abbrv}
\bibliography{hkl12}

\end{document}